\documentclass[reqno]{amsart}
\usepackage[utf8]{inputenc}

\usepackage{amsmath}
\usepackage{amssymb}
\usepackage{amsfonts}
\usepackage{amsthm}
\usepackage{multirow}
\usepackage{enumerate}
\usepackage{color}
\usepackage{dsfont}

\usepackage{setspace}
\newcommand{\R}{\mathbb{R}}
\newcommand{\C}{\mathbb{C}}
\newcommand{\f}{\rightarrow}

\newcommand{\de}{\partial}
\newcommand{\K}{K\"{a}hler}

\newcommand{\lmb}{\lambda}
\newcommand{\ov}[1]{\overline{#1}}
\newcommand{\w}[1]{\widetilde{#1}}

\newcommand{\D}{\mathcal{D}}

\newcommand{\hyp}{\operatorname{hyp}}

\newcommand{\tr}{\operatorname{tr}}

\newcommand{\arctanh}{\operatorname{arctanh}}

\newcommand{\Vol}{\operatorname{Vol}}

\newcommand{\Entp}{\operatorname{Ent_{top}}}
\newcommand{\Ent}{\operatorname{Ent_d}}
\newcommand{\Entv}{\operatorname{Ent_v}}
\newcommand{\Entvol}{\operatorname{Ent_{vol}}}

\newcommand{\ep}{{\varepsilon}}
\newcommand{\X}{{\mathcal X}}
\newcommand{\M}{{\mathcal M}}

\newcommand{\dive}{\operatorname{div}}
\newcommand{\om}{{\mathrm{hyp}}}
\newcommand{\rk}{\operatorname{rk}}
\newcommand{\di}{{\operatorname{\rho}}}

\newtheorem{thm}{Theorem}[section]
\newtheorem{prop}[thm]{Proposition}
\newtheorem{lem}[thm]{Lemma}
\newtheorem{cor}[thm]{Corollary}
\newtheorem{defn}[thm]{Definition}

\newtheorem{rmk}[thm]{Remark}

\begin{document}

\author[R. Mossa]{Roberto Mossa}

\address{Laboratoire de Math\'ematiques Jean Leray (UMR 6629) CNRS, 2 rue de la Houssini\`ere -- B.P. 92208 -- F-44322 Nantes Cedex 3, France}

\email{roberto.mossa@gmail.com}

\title[Upper and lower bounds for the first eigenv. and the vol. entropy]
{Upper and lower bounds for the first eigenvalue and the volume entropy of noncompact  K\"ahler manifolds}

\begin{abstract}
We find upper and lower bounds for the first eigenvalue and the volume entropy of a noncompact real analytic K\"ahler manifold, in terms of Calabi's diastasis function and diastatic entropy, which are sharp in the case of the complex hyperbolic space. As a corollary we obtain explicit lower bounds for the first eigenvalue of the geodesic balls of an Hermitian symmetric space of noncompact type.
\end{abstract}

\maketitle

\tableofcontents

\section{Introduction and statement of the main results}
The first eigenvalue $\lmb_1(M,g)$ of the Laplace operator of a Riemannian manifold $\left(M,\, g\right)$ is one of the most natural Riemannian invariant and its estimation is a classical problem (see e.g. \cite{cheeger, DOCARMO}).  When $M$ is noncompact with or without boundary, we call the \emph{bottom of the sprectrum} or, with abuse of language, the \emph{first eigenvalue} of $\left(M, \, g\right)$ the following Riemannian invariant
\[
\lmb_1\left(M,g\right)=\inf\left\{ \lmb_1\left(D\right) : D \in \mathcal M_0   \right\},
\]
where $\mathcal M_0$ is the set of compact domains $D$ contained in the interior of $M$ with regular boundary $\de D$  and $\lmb_1(D)$ denotes the \emph{first eigenvalue} of $D$ i.e. the \emph{smallest} $\lmb$ satisfying
\[
\Delta u = \lmb \, u
\]
for same non zero function $u$ on $M$ with $u_{|\de D}=0$.
We are interested in finding upper and lower bound of the first eigenvalue in the case of noncompact \K\ manifolds (we refer to \cite{AGL} and reference therein for the compact case).
More precisely, we consider the case of real analytic \K\ manifolds $(M,g)$ which admit a globally defined Calabi's diastasis function $\D_p: M \f \R$, the canonical \K\ potential defined by E. Calabi in his celebrated paper \cite{calabi} 
(see next section for the definition of the diastasis and its main properties).  
Our first result is the following theorem, where we give a lower bound of $\lmb_1(M,g)$ in terms of $\D_p$. 
\begin{thm}\label{thm barta}
Let $(M,g)$ be a real analytic \K\ manifold of complex dimension $n$ and let $p\in M$ be a point  for which the diastasis $\D_p:M \f \R$ is globally defined. Then
\begin{equation}\label{eq thm barta}\frac{4\, n^2}{\X\left(p\right)}\leq\lmb_1\left(M,\, g\right),
\end{equation}
where $\X \left(p\right) = \sup_q \left\| d_q \D_p \right\|^2$. If   $\X \left(p\right) = \infty$ we set $\frac{1}{\X \left(p\right)} =0$.
\end{thm}

When the manifold involved is an Hermitian symmetric space of noncompact type (HSSNT in the sequel) one obtains the following corollary of Theorem \ref{thm barta}.
\begin{thm}\label{cor hssnct}
Let $\left(\Omega, \, g_{\om}\right)$ be an $n$-dimensional irreducible HSSNT with holomorphic sectional curvature between $-4$ and $0$ and rank $r$. Then
\begin{equation}\label{eq cor HSSNT0}
\frac{n^2}{r}\leq \lmb_1\left(\Omega,\, g_{\om}\right),
\end{equation}
which is an equality if and only if $\left(\Omega,\, g_{\om}\right)$  is the complex hyperbolic space. 
 Moreover, if we denote by $B^{\Omega}_p\left(t\right)\subset\Omega$ the geodesic ball of radius $t$ and centre $p$, we have
\begin{equation}\label{eq cor HSSNT}
\frac{n^2}{r\tanh^2\left( \frac{t}{\sqrt r}\right)}\leq \lmb_1\left(B^{\Omega}_p\left(t\right),\, g_{\om}\right).
\end{equation}
\end{thm}

As a consequence of Theorem \ref{thm barta} we get a lower bound of the volume entropy $\Entv\left( M,\,  g\right)$  (see next section for details) i.e.
\begin{cor}\label{cor ent 1}\rm
Let $\left(M,\, g\right) $ be as in Theorem \ref{thm barta}. Assume moreover that $\left( M,\,  g\right)$ is complete with infinite volume. Then 
\begin{equation}\label{eq cor barta}
\frac{4\, n}{\sqrt{\X\left(p\right)}}\leq\Entv\left( M,\,  g\right),
\end{equation}
where $\Entv\left( M,\,  g\right)$ is defined by\footnote{Notice that formula \eqref{def vol ent int} make sense also in the Riemannian setting (see the discussion in Section \ref{definition varie})}
\begin{equation}\label{def vol ent int}
\Entv\left(M,\, g\right)= \inf\left\{c \in \R^+: \, \int_{M} e^{-c\, \di \left(p,\, x\right)}\, d v_g \left( x \right)  < \infty \right\}.
\end{equation}
\end{cor}

\smallskip

Inspired by the results obtained in \cite{exponential} where the author joint with A. Loi studied the concept of \emph{diastatic exponential}, defined by substituting in the usual exponential map the geodesic distance with the diastasis, it has been natural to give the following definition of \emph{diastatic entropy}, obtained by replacing in \eqref{def vol ent int} the geodesic distance with the diastasis function:

 \begin{defn}\label{defn diast ent} \rm 
Let $\left(M,\, g\right)$ be a \K\ manifold with globally defined diastasis $\D_p: M \f \R$ centred in $p \in M$, the \emph{diastatic entropy} in $p$ is defined by
\begin{equation}\label{int entr}
\begin{split}
\Ent\left(M,\, g\right)\left(p\right)=\inf\left\{c \in \R^+: \,  \int_{ M} e^{-c\, \D_p}\, d v_g < \infty\right\}.
\end{split}
\end{equation}
\end{defn}

Our last result is represented by the following theorem:
\begin{thm}\label{thm upper}
Let $\left(M,g\right)$ be a complete \K\ manifold with infinite volume and diastasis $\D_p$ globally defined for same point $p\in M$. Then 
\begin{equation}\label{eq thm upper}
\Entv \left(M,\, g\right) \leq {\Ent \left(M,\, g\right)\left(p\right)} \sqrt{ \X\left(p\right)}.
\end{equation}
and
\begin{equation}\label{eq thm upper2}
\lmb_1\left(M,\, g\right) \leq \frac{\Ent^2 \left(M,\, g\right)\left(p\right)}{4} {\, \X\left(p\right)}.
\end{equation}
\end{thm}
\smallskip

The paper is organized as follows. In the next section, after recalling the definition and the main properties of Calabi's diastasis, we show that our definition of volume entropy \eqref{def vol ent int} extends the classical one. Moreover we show that Definition \ref{defn diast ent} of the diastatic entropy, essentially, does not depend on the point $p$ and we give lower and upper bounds for the volume entropy of an HSSNT  in terms of its diastatic entropy (Corollary \ref{cor ent su giu}). In Section \ref{sect proofs} we prove Theorem \ref{thm barta} and Theorem \ref{thm upper}.  In the last section we recall some standard facts about HSSNT and Hermitian positive Jordan triple systems, needed in the proof of Theorem \ref{cor hssnct}. 
\smallskip

\noindent {\bf Acknowledgments}. The author would like to thank Professor Andrea Loi for his interest in my work and his numerous comments, Professor Gilles Carron for various stimulating discussions and Professor Sylvestre Gallot for his help and continuous encouragement.

\section{Calabi's diastasis function, volume and diastatic entropy}\label{definition varie}

\noindent {\bf Calabi's diastasis function.} 
We briefly recall the main properties of the diastatis function, defined by  Calabi in its seminal paper \cite{calabi}, the key tool of this paper.
Let $M$ be a complex manifold
endowed with a 
real analytic \K\ metric $g$.
A \K\ potential for $g$
is a real
analytic function 
$\Phi: U \f R$ 
defined in a neighborhood
of a point $p$
such that 
$\omega =\frac{i}{2}\ov\partial\partial\Phi$,
where $\omega$
is the \K\ form
associated
to $g$.
The \K\ potential $\Phi$ is not unique:
it is defined up to an addition with
the real part of a holomorphic function.
By duplicating the variables $z$ and $\ov z$ the potential $\Phi$ can be complex analytically
continued to a function 
$\w\Phi$ defined in a neighborhood of the diagonal of $ U \times \ov U$
(where $\ov U$ is the manifold
conjugated to $U$).
The {\em diastasis function} centered in $p$ is the 
\K\ potential $D_p: W \f \R$ around $p$ defined by
$$D_p(q)=\w\Phi (q, \ov q)+
\w\Phi (p, \ov p)-\w\Phi (p, \ov q)-
\w\Phi (q, \ov p).$$
We say that $D_p$ is \emph{globally defined} when $W = M$. 
The basic properties of the diastasis function (see \cite{calabi}) are the following:
\begin{enumerate}[(i)]
\item it is uniquely determined by the Kähler metric: it does not depend on the choice of the \K\ potential $\Phi$;
\item it is real valued in its domain of (real) analyticity;
\item it is symmetric in the sense that $\D_p(q)=\D_q(p)$;
\item it is equal to zero in the origin i.e. $D_p( p) = 0$;
\item if $\di_p(q)$ denotes the distance between $p\in M$ and $q$, then
\[
\D_p(q)=\di_p^2(q)+O(\di_p^4(q)) \text{ when } \di_p(q) \f 0_+ .
\]
\end{enumerate}

\smallskip

\noindent {\bf The volume entropy.} For a compact riemannian manifold $\left( X, \, g \right)$ the classical definition of volume entropy 
is the following
\begin{equation}\label{def vol ent class}
\begin{split}
\Entvol \left(X,\, g\right) =  \lim_{t \f \infty} \frac{1}{t}\log\Vol\left( \w B_p\left(t\right)\right),
\end{split}
\end{equation}
where $\Vol\left( \w B_p\left(t\right)\right)$ denotes the volume of the geodesic ball $\w B_p\left(t\right)\subset \w X$, of center $p$ and radius $t$, contained in the Riemannian universal covering $\left( \w X, \, \w g \right)$ of $\left( X, \, g \right)$. 

This notion of entropy is related with one of the main invariant for the dynamics of the geodesic flow of $\left(X,\, g\right)$, the topological entropy $\Entp \left(X,\, g\right)$ of this flow. For every compact manifold $\left(X,\, g\right)$ A. Manning in \cite{manning} proved the inequality 
$\Entvol \left(X,\, g\right) \leq \Entp \left(X,\, g\right)$, 
which is an equality when the curvature is negative. We refer the reader to the paper \cite{bcg1} (see also \cite{bcg2} and \cite{bcg3}) of  G. Besson, G. Courtois and S. Gallot  for an overview about the volume entropy and for the proof of the celebrated minimal entropy theorem.

The next lemma show that the two definition of volume entropy \eqref{def vol ent int} and \eqref{def vol ent class} coincide in the following sense
\begin{equation}\label{def vol ent class coinc}
\Entvol \left(X,\, g\right)=\Entv \left(\w X,\, \w g\right).
\end{equation}
\begin{lem}\label{classicalentropy}
Let $\left( X, \, g \right)$ be a complete n-dimensional Riemannian manifold with infinite volume. Denote by $$\underline L := \lim_{R \f + \infty}  \inf \left( \frac{1}{R} \log \left( \Vol B \left( x_0, \, R \right) \right) \right)$$ and $$\ov L :=\lim_{R \f + \infty}  \sup \left( \frac{1}{R} \log \left( \Vol B \left( x_0, \, R \right) \right) \right),$$ where $B \left( x_0, \, R \right) \subset \left( X, \, g \right)$ is the geodesic ball of centre $x_0$ and radius $R$. Then the  two limits does not depend on $x_0$ and
\begin{equation*}
\underline L \leq \Entv \left( X, \, g \right) \leq  \ov L.
\end{equation*} 
\end{lem}
\proof
The inferior limit does not depend on $x_0$, indeed set $D = d \left( x_0, \, x_1 \right)$ and $R > D$, by the triangular inequality
\[
B \left( x_0, \, R- D \right)  \subset B \left( x_1, \, R \right) \subset B \left( x_0 ,\, R+D \right)
\] 
so
\[
\lim_{R \f + \infty}  \inf \left( \frac{1}{R} \log \left( \Vol B \left( x_1, \, R \right) \right) \right) \leq
\lim_{R \f + \infty}  \inf \left( \frac{1}{R} \log \left( \Vol B \left( x_0, \, R +D \right) \right) \right) 
\]
\[
= \lim_{R' \f + \infty}  \inf \left( \frac{R'}{R'-D}\frac{1}{R'} \log \left( \Vol B \left( x_0, \, R' \right) \right) \right) 
\]
\[
\leq\lim_{R' \f + \infty}  \inf \left( \frac{1}{R'} \log \left( \Vol B \left( x_0, \, R' \right) \right) \right). 
\]
With an anologous argument one can prove the inequality in the other direction and the equality for the superior limit.

Assume $\ov L < \infty$. By the definition of limit inferior and superior, for every $\ep > 0$, there exists $R_0(\ep)$ such that, for $R \geq R_0(\ep)$,
\[
\underline L - \ep \leq \left( \frac{1}{R} \log \left( \Vol B \left( x_0, \, R \right) \right) \right) \leq \ov L + \ep
\]
equivalently
\begin{equation}\label{liminfsup}
e^{\left( \underline L - \ep \right) R}  \leq \left( \Vol B \left( x_0, \, R \right)  \right) \leq e^{\left( \ov L + \ep \right) R}.
\end{equation}
Integrating by parts we obtain
\[
I:=\int_M e^{-c\,  \di \left( x_0 ,\, x \right)} dv(x) = \int_0^\infty e^{-c\, r} \Vol_{n-1} \left( S\left( x_0 ,\, r \right) \right) dr 
\]
\[
= c \int_0^\infty e^{-c\,  r}  \Vol \left( B\left( x_0 ,\, r \right) \right) dr,
\]
where in the last equality we use that $\ov L < \infty$. Now, by \eqref{liminfsup}
\[
\int_{R_0(\ep)}^\infty e^{\left(  \underline L  -c -\ep \right)  r}\,  dr \leq \int_{R_0(\ep)}^\infty e^{-c\,  r}  \Vol \left( B\left( x_0 ,\, r \right) \right) dr \leq \int_{R_0(\ep)}^\infty e^{-\left( c - \ov L -\ep   \right)  r}\,  dr.
\]
By the second inequality one can deduce that if $c > \ov L$ then $I$ is convergent, i.e. $\ov L \geq \Entv$. On the other hand the first inequality imply that $I$ is not convergent when $c < \underline L$, that is $\Entv \geq \underline L$, as wished. 
\endproof

\smallskip

\noindent {\bf The diastatic entropy.} The following proposition shows that the Definition \ref{defn diast ent} of the diastatic entropy, essentially, does not depend on the point $p$ chosen
\begin{prop}
Let $\left(M,\, g\right)$ be a \K\ manifold and let $\D_p:M \f \R$ be globally defined for every $p\in M$ and $\sup_{p \in M}\X\left(p\right) < \infty$, then $\Ent\left( M, \, g \right)\left( p \right)$ does not depend on $p \in M$.
\end{prop}
\proof
Denoted $\ov\X = \sup_{p}\sqrt{\X\left(p\right)} = \sup_{p,q} \left\| d_q \D_p \right\|$, for every $p,\,q,\,x \in M$ we have that
\[
\left| \D_p\left(x\right) - \D_q\left(x\right) \right| = \left| \D_x\left(p\right) - \D_x\left(q\right) \right| \leq \ov\X \, \di\left( p, \, q \right).
\]
So
\[
e^{-c \, \ov\X\, \di\left( p, \, q \right)} \int e^{-c \, \D_p \left( x \right)} \leq  \int e^{-c \, \D_q \left( x \right)} \leq e^{c \, \ov\X\, \di\left( p, \, q \right)} \int e^{-c \, \D_p \left( x \right)}
\]
Therefore $\int e^{-c \, \D_q \left( x \right)} < + \infty$ if and only if $\int e^{-c \, \D_p \left( x \right)} < + \infty$.

\endproof

In \cite{mossa1}  the author show that $\Ent\left(M,\, g\right)(p)$ and  the balanced condition on $g$ are deeply linked (the notion of balancedness has been defined by S. K. Donaldson \cite{donaldson} for the compact case and then extended to the noncompact case by C. Arezzo and A. Loi \cite{arezzoloi}). In the same paper one can also find the computation of the diastatic entropy for every  homogeneous domain in terms of Piatetskii-Shapiro constants \cite[(8)]{mossa1}. In particular, for any irreducible HSSNT,  we have
\begin{equation}\label{entd hssnct}
\begin{split}
\Ent\left(\Omega,\, g_{\om}\right)=\gamma - 1, 
\end{split}
\end{equation}
where $\gamma$ is the genus of $\Omega$. If $g_B$ denotes the Bergman metric on $\Omega$, the genus is defined by $g_B = \gamma \, g_\om$. This formula should be compared with the expression of the volume entropy (see \cite{mossa})
\begin{equation}\label{entv hssnct}
\operatorname{Ent_v}\left(\Omega,\, g_{\om}\right)=2\,\sqrt{\sum _{j=1}^{r}
 \left( b+1+a \left( r-j \right)  \right) ^{2}},
\end{equation}
 in terms of the invariants $r,$ $a$ and $b$ associated to $\Omega$ (see \cite[Table 1]{mossa} for a description of these invariants). Recalling that $\gamma=b+ 2 + a \left( r-1\right)$ (see \cite{mossa1}) we obtain the following result, (which is in accordance with \eqref{eq thm upper}, since, by \eqref{eq chi} below, $\X(p )=4\, r$)
\begin{cor}\label{cor ent su giu}
Let $\left(\Omega, \, g_{\om}\right)$ be an $n$-dimensional irreducible HSSNT with holomorphic sectional curvature between $-4$ and $0$ and rank $r$. We have
\begin{equation}\label{ineqentventdhssnct}
2\, {\Ent \left(\Omega,\, g_{\om}\right)} \leq \Entv \left(\Omega,\, g_{\om}\right) \leq 2\,\sqrt r \, {\Ent  \left(\Omega,\, g_{\om}\right)},
\end{equation}
where the equalities are attained only for $r=1$ or $a=0$ i.e. when $\Omega$ is the complex hyperbolic disc.
\end{cor}

\section{Proof of Theorems \ref{thm barta} and \ref{thm upper}}\label{sect proofs}

The main ingredient in the proof of Theorems \ref{thm barta} is the following lemma. 
\begin{lem}\label{lemma barta}(Barta's Lemma, \cite{barta})
Let $\phi:M \f \R$ be a positive function such that $\Delta \phi \geq \lmb\, \phi$, then 
\[
\lmb_1\left(M\right)\geq \lmb.
\]  
\end{lem}
\proof For the sake of completeness we will give a proof of this lemma. 
Given  $f \in \mathcal C^1_0(M)$, we can write it as $f=u\,\phi$, where $u \in  \mathcal C^1_0(M)$. We have
\[
\|df\|^2= \phi^2 \,  \| du \|^2+u^2  \, \| d \phi \|^2+\frac{1}{2} \, \langle d \phi^2, \,  d u^2 \rangle.
\]
As 
\[
\dive \left(u^2 \,  d \phi^2\right)=-\tr \nabla \left(u^2 \,  d \phi^2\right) = u^2 \,  \Delta \phi^2 - \langle du^2, \, d\phi^2 \rangle
\]
we get
\begin{equation*}\|df\|^2=\phi^2 \,  \| du \|^2+u^2 \left( \| d \phi \|^2 + \frac{1}{2} \,  \Delta \phi^2 \right) - \frac{1}{2}\,\dive \left(u^2 \,  d \phi^2\right)
\end{equation*}
\begin{equation*}
=\phi^2 \,  \| du \|^2+u^2 \,  \phi  \, \Delta \phi - \frac{1}{2}\,\dive \left(u^2  \, d \phi^2\right).
\end{equation*}
where we used that $\frac{1}{2}\, \Delta \phi^2= \phi\, \Delta \phi - \|d \phi\|^2$. Since $u^2 \,  d \phi^2$ and $df$ have compact support we have $\int_M \dive (u^2 \, d \phi^2)=0$ and $\int_M \|df\|^2<\infty$, therefore
\begin{equation*}
\int_M \|df\|^2= \int_M \left(\phi^2 \, \|d u\|^2+ u^2  \, \phi\, \Delta \phi\right)
\end{equation*}
\begin{equation*}
\geq\int_M u^2 \,  \phi\, \Delta \phi=\lmb  \, \int_M u^2 \,  \phi^2 = \lmb \,  \int_M f^2.
\end{equation*}
Then we have shown that
\[
\lmb \leq \frac{\int_M \|df\|^2}{\int_M |f|^2},
\]
for every $f \in \mathcal C^1_0(M)$. The inequality $\lmb_1\left(M\right)\geq \lmb.$ is a consequence of the equality
\begin{equation}\label{lemma00}
\lmb_1\left(M,g \right)=\inf_{f\in \mathcal
H^1_0(M)} \frac{\int_M\|df\|^2}{\int_M|f|^2},
\end{equation}
were $\mathcal H^1_0\left(M, \, g\right)$ denotes the completion of the space of the $\mathcal C^1$-differentiable function with compact support contained in the interior of $M$ with respect the norm $\|f\|_{1,2} = \sqrt{ \int_M |f|^2 + \int_M \left| df \right|^2 }$ (see \cite[Lemma D.II.3]{berger} for a proof).
\hfill $\square$

\bigskip

\noindent {\bf Proof of Theorem \ref{thm barta}.}
Set $\phi:=e^{-c\, \D_p}$. We have: 
\[
\Delta \phi = - \tr (\nabla d \phi)= c  \left(\tr \left(\nabla d \D_p\right) - c\,\|d \D_p\|^2\right)\phi,
\]
by the identity $\nabla d \D_p (v,v) + \nabla d \D_p (Jv,Jv)= 4 g(v,v)$ we get
\[
\Delta \phi \geq c  \left(4\,n -  c\, \X(p) \right)\phi,
\]
then by Barta's Lemma we conclude
\[
\lmb_1 \geq \max_{c\in \R}\left[ c \left( 4\,n - c\, \X(q) \right) \right] =\frac{4  \, n^2}{\X(p)}
\]
and this ends the proof of Theorem \ref{thm barta}.
\hfill $\square$

\bigskip

\noindent {\bf Proof of Theorem \ref{thm upper}.}
Since $\D_p\left(p\right) =0$, for every $x \in M$ we have 
\[
\D_p\left(x\right) = \D_p\left( x\right) - \D_p\left(p\right) \leq \sup_{z \in M} \left\| d_ z \D_p \right\| \di_p \left( x \right) \leq \sqrt{\X (p)}\ \di_p\left( x \right),
\]
so
\[
\int_M e^{-c\, \X \left( p \right) \, \di_p\left( x \right)}  \leq \int_M e^{-c\, \D_p \left( x \right)}.
\]
We deduce that if $c \, \sqrt {\X \left( p \right)} \leq \Entv \left(M,\, g\right) $ then $c \leq \Ent \left(M,\, g\right)$. In particular, taking $c=\frac {\Entv \left(M,\, g\right)}{\sqrt{ \X \left( p \right)}}$, we obtain \eqref{eq thm upper}.

Now set $f(z):=e^{-c\,\di_x(z)}$, is not hard to prove that $f \in  \mathcal H^1_0(M)$ for every $c > \frac{\Entv (M, \, g)}{2}$. Substituting in \eqref{lemma00} we get
\[
\lmb_1\leq \frac{\int_{M}\|df\|^2}{\int_{M}|f|^2}= c^2
\]
where in the last equality we use that $\left\| d \rho_x \right\|^2=1$. As $c$ approach to $ \frac{\Entv (M, \, g)}{2}$ we obtain the following inequality 
\begin{equation}\label{eq entv lmb}
\lmb_1 \left(M, g\right) \leq \frac{\Entv^2 \left(M,\, g\right)}{4}
\end{equation}
(which is an equality when $\left(M,\, g\right)=\left(\Omega,\, g_{\om}\right)$, see \cite[Appendix C]{bcg0} for a proof). 
Finally combining \eqref{eq thm upper} with \eqref{eq entv lmb} we obtain \eqref{eq thm upper2}.

\begin{flushright} $\square$ \end{flushright}

\begin{rmk}\label{prove inequality eq entv lmb}\rm
Observe that by the same argument used to prove \eqref{eq entv lmb} one can obtain an alternative proof of inequality \eqref{eq thm upper2} by setting $f\left( z \right) := e^{-c \, \D_p(z)}$.
\end{rmk}

\section{Hermitian positive triple system and proof of Theorem \ref{cor hssnct}}\label{app}

We refer the reader  to \cite{roos} (see also \cite{loos}) for more details on Hermitian symmetric spaces of noncompact type and Hermitian positive Jordan triple systems (from now on HPJTS).

\smallskip

\noindent {\bf  Definitions and notations.}
An Hermitian Jordan triple system is a  pair $\left({\mathcal M},
\{ ,  ,\}\right)$, where ${\mathcal M}$ is a complex vector space and $\{
,  ,\}$ is a map
\[
\{ ,  ,\}:{\mathcal M}\times {\mathcal M}\times {\mathcal M} \rightarrow {\mathcal M}
\]
\[
\left(u, \,  v, \,  w\right)\mapsto \{ \, u, \,  v, \,  w\}
\]
which is ${\C}$-bilinear and symmetric in $u$ and $w$, ${\C}$-antilinear in $v$ and such that the following \emph{ Jordan identity} holds:
\[
\{x,  \, y, \,  \{u, \,  v,  \, w\}\}-\{u,  \, v, \,  \{x,  \, y, \,  w\}\}= \{\{x, \,  y, \,  u\}, \,  v,
 \, w\}-\{u, \,  \{v,  \, x,  \, y\},  \, w\}.
\]
For $x, \, y, \, z \in \M$ consider  the following operators
\[
T\left(x, \, y\right)z =\left\{  x, \, y, \, z\right\} 
\]
\[
Q\left(x, \, z\right) \, y =\left\{  x, \, y, \, z\right\}  
\]
\[
Q\left(x, \, x\right) =2\,Q\left(x\right)\label{D3}\\
\]
\[
B\left(x, \, y\right) =\operatorname{id}_{\mathcal M}-T\left(x, \, y\right)+Q\left(x\right)Q\left(y\right). \label{D4}
\]
The operators $B\left(x, \, y\right)$ and $T\left(x, \, y\right)$ are $\mathbb{C}$-linear and the operator
$Q\left(x\right)$ is $\mathbb{C}$-antilinear. $B\left(x, \, y\right)$ is called \emph {Bergman operator}.
An Hermitian Jordan triple system is called  \emph {positive} if the Hermitian form
\[
\left(  u\mid v\right)  =\tr T\left(u, \, v\right) \label{D5}
\]
is positive definite. An element $c \in \M$ is called \emph {tripotent} if
$\{c, \, c, \, c\}=2 \, c$. Two tripotents $c_1$ and $c_2$ are called \emph {(strongly)
orthogonal} if $T\left(c_1,  \, c_2\right)=0$.

\smallskip

\noindent{\bf HSSNT associated to HPJTS.}
M. Koecher (\cite{Koecher1}, \cite{Koecher2}) discovered that to every HPJTS
$\left(\M, \{ ,  ,\}\right)$ one can associate an HSSNT, in its realization as a bounded symmetric domain $\Omega$
centered at the origin $0\in \M$. The domain $\Omega$ is defined as the connected component containing the origin of   the set of all $u\in {\M}$ such that $B\left(u,  \, u\right)$ is positive definite with respect to the Hermitian form $\left(  u\mid v\right)  =\tr T\left(u, \, v\right) \label{D5}$. \emph{We will always consider such a domain in its (unique up to linear isomorphism) circled realization.} Suppose that $\M$ is \emph{simple} (i.e. $\Omega$ is irreducible).
 The \emph{flat} form $\omega_0$ is defined by
\[
\omega_0= -\frac{i}{2 \, \gamma} \,  \partial \ov \partial \left(  z\mid z\right),
\]
where $\gamma$ is the genus of $\Omega$.
If $\left(  z_{1},\ldots,z_{n}\right)$ are orthonormal
coordinates for the Hermitian product $\left(  u\mid v\right)  $, then%
\[
\omega_{0}=\frac{i}{2 \, \gamma}\sum_{m=1}^{d}{d}z_{m}\wedge
{d}\overline{z}_{m}.
\]
The reproducing kernel ${K_\Omega}$ of $\Omega$, with respect $\omega_0$ is given by
\begin{equation}\label{KOB}
\left(K_\Omega\left(z,  \, \ov z\right)\right)^{-1}=C \,  \det B\left(z, \, z\right),
\end{equation}
where $C=\int_\Omega \frac{\omega_0^n}{n!}$.
When $\Omega$ is irreducible,
\[
\omega_{\om}= -\frac{i}{2 \, \gamma}\, \partial \ov \partial \log\det B.
\]
is the \emph{hyperbolic} form on $\Omega$, with the associated hyperbolic metric $g_{\om}$ (whose holomorphic sectional curvature is between $-4$ and $0$). The HSSNT associated to $\M$ is $(\Omega,  \, g_\om)$. Moreover the hyperbolic form is related to the flat form by
\begin{equation}\label{flathyp}
\omega_\om(z)(u,  \, v) = \omega_0(B(z,  \, z)^{-1}u, \,  v).
\end{equation}

The HPJTS $\left({\M}, \{ ,  ,\}\right)$ can be recovered by its
associated HSSNT $\Omega$ by defining ${\M}=T_0 \Omega$ (the tangent space to the origin of $\Omega$) and
\begin{equation*}\label{trcurv}
\{u,  \, v, \,  w\}=-\frac{1}{2} \, \left(R_0\left(u,  \, v\right) \, w+J_0 \, R_0\left(u, \,  J_0\,v\right)w\right),
\end{equation*}
where $R_0$ (resp. $J_0$) is the curvature tensor of the Bergman metric (resp. the complex structure) of $\Omega$ evaluated at the origin.
For more informations on the correspondence between
HPJTS and HSSNT we  refer also  to p. 85 in Satake's book \cite{satake}.

\smallskip

\noindent {\bf  Spectral decomposition and polar coordinates.}
Let $\M$ be a HPJTS. Each element $z\in \M$ has a unique \emph{spectral decomposition}
\[
z=\lambda_{1} \, c_{1}+\cdots+\lambda_{s} \, c_{s}\qquad\left(0<\lambda_{1}<\cdots
<\lambda_{s}\right), \label{D7}
\]
where $\left(c_{1},\ldots, \, c_{s}\right)$ is a sequence of pairwise
orthogonal tripotents and the $\lambda_j$ are real number called eigenvalues of $z$. The integer $s= \rk( z)$ is called \emph{rank} of $z$.  For every $z \in \M$ let $\max\{z\}$ denote the largest eigenvalue of $z$, then $\max\{\cdot \}$ is a norm on $\M$ called the \emph{spectral norm}. The HSSNT $\Omega\subset\M$ associated to $\M$
is the open unit ball in $\mathcal M$ centered at the origin (with respect the spectral norm),
i.e.,
\begin{equation*}\label{Mball}
{\Omega}=\{z=\sum_{j=1}^s\lambda_j \, c_j \ |\ \max\{z\}= \max_j\{ \lmb _j\}<1\}.
\end{equation*}
The rank of $\M$ is $\rk( \M) = \max \{\rk( z)\, |\, z \in \M    \}$, moreover $\rk (\M) =\rk( \Omega) = :r$.
The elements $z$ such that $\operatorname{rk}\left( z \right) = r$ are called \emph{regular}. If $z\in {\M}$ is regular, with spectral decomposition%
\begin{equation*}
z=\lambda_{1} \, e_{1}+\cdots+\lambda_{r} \, e_{r}\qquad(\lambda_{1}>\cdots
>\lambda_{r}>0), \label{D8}%
\end{equation*}
then $\left(  e_{1},\ldots, \, e_{r}\right)  $ is a \emph{(Jordan) frame of }${\M}$.

The set $\mathcal{F}$ of frames (also called F\"{u}rstenberg-Satake boundary
of $\Omega$)
is a compact
manifold. The map $F: \left\{  \lambda_{1}>\cdots>\lambda_{r}>0\right\}  \times\mathcal{F}
\rightarrow {\M}_{\mathrm{reg}}\nonumber$ defined by:
\begin{equation}\label{polar coord1}
\left(  \left(  \lambda_{1},\ldots, \, \lambda_{r}\right)  ,\left(  c_{1}%
,\ldots, \, c_{r}\right)  \right)  \mapsto\sum\lambda_{j} \, c_{j}%
\end{equation}
is a diffeomorphism onto the open dense set ${\M}_{\mathrm{reg}}$ of regular elements of
${\M}$, moreover its restriction%
\begin{equation}\label{polar coord}
\left\{  1>\lambda_{1}>\cdots>\lambda_{r}>0\right\}  \times\mathcal{F}%
\rightarrow\Omega_{\mathrm{reg}}%
\end{equation}
is a diffeomorphism onto the (open dense) set $\Omega_{\mathrm{reg}}$ of regular elements
of $\Omega$. This map plays the same role as polar coordinates in rank one.

\smallskip

\noindent {\bf Proof of Theorem \ref{cor hssnct}}\label{proof cor hssnct}
By using the rotational symmetries of $\Omega\subset\M$
one can show that the diastasis function at the origin
$\D^{\hyp}_0:M\rightarrow\R$, associated to $g_{\om}$, is globally defined and  reads as
\[
\D^{\hyp}_0\left(z\right)=\frac{1}{\gamma}\log \left( C \,  K_\Omega\left(z, \ov z\right)\right),
\]
where $C=\int_\Omega \frac{\omega_0^n}{n!}$,
(see \cite{loi} for a proof and  further results on Calabi's function
for HSSNT).

Let $z \in \Omega$ be a regular point and let $z=\sum_{j=1}^r \lmb_j  \, c_j$ be its expression in polar coordinates. We have that  (see \cite{roos} for a proof)
\begin{equation}\label{BZZ}
B\left(z, \, z\right)c_j=\left(1-\lambda_j^2\right)^2c_j, \quad j=1,\dots , \, r,
\end{equation}
\[
\det B\left(z, \, z\right) = \prod_{j=1}^r\left(1-\lmb_j^2\right)^\gamma,
\]
Thus (\ref{KOB}) yields the expression of the diastasis with respect the coordinates \eqref{polar coord},
\begin{equation*}
\D^{\hyp}_0\left(z\right)=-\frac{1}{\gamma}  \, \log \det B\left(z, \, z\right)=- \log \prod_{j=1}^r\left(1-\lmb_j^2\right)
\end{equation*}
so
\begin{equation}\label{doz}
d_z\D_0^\om= \sum_{j=1}^r \frac{2\,\lmb_j}{1-\lmb_j^2}\,d\lmb_j.
\end{equation}
By \eqref{flathyp},  \eqref{BZZ} and \eqref{doz}, we see that $g_\om \left( \frac {\de}{\lmb_j}, \, \frac {\de}{\lmb_k} \right)= \frac {\delta_{jk}} {\left(  1 - \lmb_j^2  \right)}  $, we conclude that
\begin{equation}\label{ineqq norm}
\|d_z\D_0^\om\|^2= \sum_{j=1}^r \frac{4 \, \lmb^2_j}{\left(1-\lmb_j^2\right)^2} \, \|d\lmb_j\|^2=4 \, \sum_{j=1}^r \lmb_j^2.
\end{equation}

In polar coordinates of $\Omega$, the expression of the distance from the origin, with respect $g_{\om}$, is 
$\di_\om(0, \, z)=\sqrt{\sum^r_{j=1} \arctanh^2\left(\lmb_j\right) }$  (see \cite[(8)]{mossa}).
Assume now that $z \in B_0^{\Omega}(t)$, i.e. that $\sqrt{\sum_{j=1}^r \arctanh^2(\lmb_j)} < t$. 
By the concavity of the function $\arctanh^2\left(\sqrt{ |\cdot|}\right)$ we get:
\[
\sqrt r  \, \arctanh \left( \frac{\sqrt{\sum_{j=1}^r \lmb_j^2}}{\sqrt r} \right)\leq\sqrt{\sum_{j=1}^r \arctanh^2 (\lmb_j)}<t,
\]
thus
\[
\sum_{j=1}^r \lmb_j^2 \leq r  \tanh^2 \left( \frac{t}{\sqrt r}\right).
\]
Substituting the previous inequality in \eqref{ineqq norm} we get
\begin{equation}\label{eq chi}
\X(0) = \sup_{z \in B_0^{\Omega}(t)} \|d_z\D_0\|^2 \leq 4\,r  \tanh^2 \left( \frac{t}{\sqrt r}\right).
\end{equation}
Finally Theorem \ref{thm barta} yields, 
\begin{equation*}
\frac{n^2}{r\tanh^2\left( \frac{t}{\sqrt r}\right)}\leq \lmb_1\left(B^{\Omega}_p(t), \,  g_\om\right),
\end{equation*}
for any $p\in \Omega$. This shows the validity of \eqref{eq cor HSSNT}.

As $t$ tend to infinity we see that $\X(0)=\frac{n^2}{r}$, so substituting this equality in \eqref{eq thm barta} we obtain \eqref{eq cor HSSNT0}. Recalling that $n=r \left( b+1+\frac{a}{2}\left(r-1\right) \right)$, by \eqref{entv hssnct}, one can prove that
\[
\operatorname{Ent_v}\left(\Omega,\, g_{\om}\right)= \frac{2\,n}{\sqrt r}\,\sqrt{1 + \frac{a^2 \, r^2}{12\, n^2} \left( r^2 -1 \right) },
\]
therefore the equality in \eqref{eq cor HSSNT0} is attained if and only if $r=1$ or $a=0$ i.e. when $\Omega$ is the complex hyperbolic disc and this concludes the proof of Theorem \ref{cor hssnct}.
\hfill $\square$

\begin{rmk}\label{rmk cheeger}\rm
Consider $B^e_y(t)=\{x \in \C^n : \left\|z\right\|<t\}$ the ball of $\C^n$ of radius $t$ with the induced euclidean metric $g_e$. The diastasis $\D_0^e$ centered at the origin is given by $\D^e_0(z)=\|z\|^2$. By Theorem \ref{thm barta} we get
\[
\frac{n^2}{t^2} \leq \lmb_1\left(B^e_y(t), \, g_e\right),
\]
which is in accordance with the celebrated J. Cheeger inequality 
$\frac{h^2 \left(M, \, g\right)}{4}\leq\lmb_1\left(M, \, g\right)$,
where $h\left(M, \, g\right)$ is the so called Cheeger constant, defined as the minimum of $\frac{\Vol_{n-1}\left(\de D\right)}{\Vol\left(D\right)}$ among the compact domain $D$ of a Riemannian manifold $\left( M, \, g \right)$ and the isoperimetric inequality give us
$
h\left(B_y^e(t), \,  g_e\right)
=\frac{\Vol_{n-1}(\de B^e_y\left(t\right), \, g_e)}{\Vol \left(B^e_y\left(t\right), \, g_e\right)}=\frac{2\,n}{t}$.

It is also worth pointing out that X. Wang in a recent paper \cite{XWang} showed that for any HSSNT  
$\lmb_1\left(\Omega,\, g_{\om}\right) = \frac{h^2 \left(\Omega,\, g_{\om}\right) }{4}.$
Nevertheless, in general,  the computation of the Cheeger constant is a difficult matter and we will study the link between the diastasis function and the Cheeger constant in a forthcoming paper.
\end{rmk}

\end{document}